\documentclass[envcountsect]{svproc}

\usepackage{url}

\usepackage[utf8]{inputenc}
\usepackage[T1]{fontenc}
\usepackage[english]{babel}
\usepackage[autostyle=true]{csquotes}

\usepackage{amsthm}
\usepackage{latexsym, amssymb, amsmath}
\makeatletter
\renewcommand*\env@matrix[1][*\c@MaxMatrixCols c]{%
	\hskip -\arraycolsep
	\let\@ifnextchar\new@ifnextchar
	\array{#1}}
\makeatother
\usepackage{url}
\usepackage{float}
\usepackage{caption}
\usepackage{subcaption}
\usepackage{color}  
\usepackage{hyperref}
\hypersetup{
	colorlinks=true,
	linkcolor=black,
	filecolor=magenta,      
	urlcolor=cyan,
	linktoc=all, 
}
\usepackage[ruled,vlined]{algorithm2e}
\usepackage{tikz-cd}
\usepackage{graphicx}
\graphicspath{ {./images/} }
\usepackage{amscd}

\newcommand{\F}{\mathbb{F}}

\spnewtheorem{thm}{Theorem}[section]{\bfseries}{\itshape}

\newtheorem{defi}[thm]{Definition}
\newtheorem{rem}[thm]{Remark}

\newtheorem{ex}[thm]{Example}

\begin{document}
\mainmatter              
\title{Biderivations of low-dimensional \\ Leibniz algebras}
\titlerunning{Biderivations of low-dimensional Leibniz algebras}

\author{Manuel Mancini\thanks{Supported by University of Palermo and by the “National Group for Algebraic and Geometric Structures, and their Applications” (GNSAGA – INdAM).}}  
\authorrunning{Manuel Mancini} 
\tocauthor{Manuel Mancini}

\institute{Dipartimento di Matematica e Informatica\\ Universit\`a degli Studi di Palermo, Via Archirafi 34, 90123 Palermo, Italy\\
\email{manuel.mancini@unipa.it}, ORCID: 0000-0003-2142-6193.}

\maketitle              

\begin{abstract}
In this paper we give a complete classification of the Leibniz algebras of biderivations of right Leibniz algebras of dimension up to three over a field $\F$, with $\operatorname{char}(\F) \neq 2$. We describe the main properties of such class of Leibniz algebras and we also compute the biderivations of the four-dimensional \emph{Dieudonné Leibniz algebra} $\mathfrak{d_1}$. Eventually we give an algorithm for finding derivations and anti-derivations of a Leibniz algebra as pair of matrices with respect to a fixed basis.


\end{abstract}

\section*{Introduction}

Leibniz algebras were first introduced by J.-L.\ Loday in \cite{loday1993version} as a non-antisymmetric version of Lie algebras, and many results of Lie algebras were also established in the frame of Leibniz algebras. Earlier, such algebraic structures had been considered by A.\ Blokh, who called them D-algebras \cite{blokhLie}, for their strict connection with the derivations. Leibniz algebras play a significant role in different areas of mathematics and physics.

In \cite{low_dim} derivations of low-dimensional Leibniz algebras have been classified and studied. After two short preliminary sections, in this paper we aim to study and classify the Leibniz algebras of biderivations of low-dimensional Leibniz algebras. Independently of its intrinsic interest, biderivations find concrete applications in representation theory (cf. \cite{nonabext} and \cite{Casas}), (sub-) Riemannian geometry and control theory (see \cite{BiggsNagy}, and the bibliography therein), just to give two instances.

The first section is devoted to some background material on Leibniz algebras which will be useful for the rest of the paper. We address the reader to \cite{ayupov2019leibniz} and \cite{erdmann2006introduction} for more details.

In Section 2 we give the definitions of \emph{derivation}, \emph{anti-derivation} and \emph{biderivations} for a right Leibniz algebra, which have been first introduced by J.-L.\ Loday in \cite{loday1993version}. We also recall the main properties of the right and left adjoint maps and, as an example, we show how to compute the Leibniz algebra of biderivations of the four-dimensional \emph{Dieudonné Leibniz algebra} $\mathfrak{d_1}$ (cf. \cite{LM2022}). 

In the last section we classify the Leibniz algebra of biderivations of two and three-dimensional right Leibniz algebras over a field $\F$, with $\operatorname{char}(\F) \neq 2$. We use the classification of low-dimensional complex Leibniz algebras and their derivations (cf.\ \cite{ayupov2019leibniz} and \cite{low_dim}) and we give an algorithm for finding derivations and anti-derivations of a Leibniz algebras as pair of matrices with respect to a fixed basis. 

\section{Leibniz algebras}

We assume that $\F$ is a field with $\operatorname{char}(\F)\neq2$. For the general theory we refer to \cite{ayupov2019leibniz}.

\begin{defi}
	A \emph{right Leibniz algebra} over $\mathbb{F}$ is a vector space $L$ over $\mathbb{F}$ endowed of a bilinear map (called $commutator$ or $bracket$) $\left[-,-\right]\colon L\times L \rightarrow L$ which satisfies the \emph{right Leibniz identity}
	$$
	\left[\left[x,y\right],z\right]=\left[[x,z],y\right]+\left[x,\left[y,z\right]\right],\;\; \forall x,y,z\in L.
	$$
	
\end{defi}
In the same way we can define a left Leibniz algebra, using the left Leibniz identity
$$
\left[x,\left[y,z\right]\right]=\left[\left[x,y\right],z\right]+\left[y,\left[x,z\right]\right], \;\; \forall x,y,z\in L.
$$

A Leibniz algebra that is both left and right is called \emph{symmetric Leibniz algebra}. From now on we assume that $\dim_\F L<\infty$.

Clearly every Lie algebra is a Leibniz algebra and every Leibniz algebra with skew-symmetric commutator is a Lie algebra. Thus it is defined an adjunction (see \cite{mac2013categories}) between the category $\textbf{LieAlg}_{\F}$ of the Lie algebras over $\F$ and the category $\textbf{LeibAlg}_{\F}$ of the Leibniz algebras over $\F$. The left adjoint of the inclusion $i\colon\textbf{LieAlg}_\F\rightarrow \textbf{LeibAlg}_\F$ is the functor $\pi\colon\textbf{LeibAlg}_\F \rightarrow \textbf{LieAlg}_\F$ that associates with every Leibniz algebra $L$ the quotient $L/\operatorname{Leib}(L)$, where $\operatorname{Leib}(L)=\langle \left[x,x\right] |\,x\in L \rangle$ is called the \emph{Leibniz kernel} of $L$. We observe that $\operatorname{Leib}(L)$ is the smallest bilateral ideal of $L$ such that $L/\operatorname{Leib}(L)$ is a Lie algebra. Moreover $\operatorname{Leib}(L)$ is an abelian algebra. 

We define the left and the right center of a Leibniz algebra
$$
\operatorname{Z}_l(L)=\left\{x\in L \,|\,\left[x,L\right]=0\right\}, \;\; \operatorname{Z}_r(L)=\left\{x\in L \,|\,\left[L,x\right]=0\right\}
$$
and we observe that they coincide when $L$ is a Lie algebra. The \emph{center} of $L$ is $\operatorname{Z}(L)=\operatorname{Z}_l(L)\cap \operatorname{Z}_r(L)$. In the case of symmetric Leibniz algebras, the left center and the right center are bilateral ideals, but in general $\operatorname{Z}_r(L)$ is an ideal of the right Leibniz algebra $L$, meanwhile the left center is not even a subalgebra.

Finally we recall the definitions of nilpotent and solvable Leibniz algebras.

\begin{defi}
	Let $L$ be a right Leibniz algebra over $\F$ and let
	$$
	L^{(0)}=L,\,\,L^{(k+1)}=[L^{(k)},L],\;\;\forall k\geq0,
	$$ be the \emph{lower central series of $L$}. $L$ is \emph{$n-$step nilpotent} if $L^{(n-1)}\neq0$ and $L^{(n)}=0$.
\end{defi}

\begin{defi}
	Let $L$ be a right Leibniz algebra over $\F$ and let
	$$
	L^{0}=L,\,\,L^{k+1}=[L^{k},L^k],\;\;\forall k\geq0,
	$$ be the \emph{derived series of $L$}. $L$ is \emph{$n-$step solvable} if $L^{n-1}\neq0$ and $L^{n}=0$.
\end{defi}

\section{Derivations, anti-derivations and biderivations}

In this section we recall the definitions of \emph{derivation}, \emph{anti-derivation} and \emph{biderivation} for right Leibniz algebras and we show an example of computation of the biderivations algebra. 

The definition of \emph{derivation} is the same as in the case of Lie algebras.
\begin{defi}
	Let $L$ be a Leibniz algebra over $\F$. A \emph{derivation} of $L$ is a linear map $d\colon L\rightarrow L$ such that 
	$$
	d(\left[x,y\right])=\left[d(x),y\right]+\left[x,d(y)\right],\;\; \forall x,y\in L.
	$$
\end{defi}

\begin{rem}
	Fixed a basis $\lbrace e_1,\ldots,e_n \rbrace$ of $L$, where $n=\dim_{\F}L$, a linear map $d\colon L \rightarrow L$ is a derivation if and only if
	$$
	d(\left[e_i,e_j\right])=\left[d(e_i),e_j\right]+\left[e_i,d(e_j)\right], \,\,\forall i,j=1,\ldots,n.
	$$
\end{rem}

The right multiplications are particular derivations called \emph{inner derivations} and an equivalent way to define a right Leibniz algebra is to saying that the (right) adjoint map $\operatorname{ad}_x=\left[-,x\right]$ is a derivation, for every $x\in L$. Meanwhile, for a right Leibniz algebra, the left adjoint maps are not derivations in general.\\ With the usual bracket \hbox{$\left[d_1,d_2\right]=d_1\circ d_2 - d_2\circ d_1$}, the set $\operatorname{Der}(L)$ is a Lie algebra and the set  $\operatorname{Inn}(L)$ of all inner derivations of $L$ is an ideal of $\operatorname{Der}(L)$. Furthermore, $\operatorname{Aut}(L)$ is a Lie group and the associated Lie algebra is $\operatorname{Der}(L)$.\\

The definitions of \emph{anti-derivation} and \emph{biderivation} for a Leibniz algebra were first given by J.-L.\ Loday in \cite{loday1993version}.

\begin{defi}
	An \emph{anti-derivation} of a right Leibniz algebra $L$ is a linear map $D\colon L \rightarrow L$ such that
	$$
	D([x,y])=[D(x),y]-[D(y),x], \; \; \forall x,y \in L.
	$$
\end{defi}
\noindent
For a left Leibniz algebra we have to ask that
$$
D([x,y])=[x,D(y)]-[y,D(x)], \; \; \forall x,y \in L.
$$
We observe that in the case of Lie algebras, there is no difference between a derivation and an anti-derivation. Moreover one can check that, for every $x \in L$, the left adjoint map $\operatorname{Ad}_x=[x,-]$ defines an anti-derivation.

\begin{rem}
	The set of anti-derivations of a Leibniz algebra $L$ has a $\operatorname{Der}(L)-$module structure with the action
	$$
	d \cdot D := [D,d]= D \circ d - d \circ D,
	$$
	for every $d \in \operatorname{Der}(L)$ and for every anti-derivation $D$.
\end{rem}

\begin{rem}
	Let $D\colon L \rightarrow L$ be an anti-derivation. Then, for every $x \in L$, we have
	$$
	D([x,x])=[D(x),x]-[D(x),x]=0,
	$$
	thus $D(\operatorname{Leib}(L))=0$.
\end{rem}

\begin{defi}
	Let $L$ be a right Leibniz algebra. A \emph{biderivation} of $L$ is a pair
	$$
	(d,D)
	$$
	where $d$ is a derivation and $D$ is an anti-derivation, such that
	\begin{equation}\label{1}
		[x,d(y)]=[x,D(y)], \; \; \forall x,y \in L.
	\end{equation}
\end{defi}

The set of all biderivations of $L$, denoted by $\operatorname{Bider}(L)$, has a Leibniz algebra structure with the bracket
$$
[(d,D),(d',D')]=(d \circ d' - d' \circ d, D \circ d' - d' \circ D), \; \; \forall (d,D),(d,D') \in \operatorname{Bider}(L),
$$
and it is possible to define a Leibniz algebra morphism
$$
L \rightarrow \operatorname{Bider}(L)
$$
by
$$
x \mapsto (-\operatorname{ad}_x, \operatorname{Ad}_x), \; \; \forall x \in L.
$$
The pair $(-\operatorname{ad}_x, \operatorname{Ad}_x)$ is called \emph{inner biderivation} of $L$ and the set of all inner biderivations forms a Leibniz subalgebra of $\operatorname{Bider}(L)$.

Now we give an example of computation of the biderivations of a Leibniz algebra.

\begin{ex}
	Let $\mathfrak{d}_1$ be the four-dimensional \emph{Dieudonné Leibniz algebra} (see \cite{LM2022} for more details), i.e.\ $\mathfrak{d}_1$ has basis $\lbrace e_1,e_2,e_3,z \rbrace$ and non-zero commutators
	$$
	[e_1,e_3]=[e_2,e_3]=-[e_3,e_1]=[e_3,e_2]=z.
	$$
	We want to find the Leibniz algebra $\operatorname{Bider}(\mathfrak{d}_1)$ of biderivations of $\mathfrak{d}_1$. The derivations of the Dieudonné Leibniz algebra $\mathfrak{d_n}$ over a field $\F$, with $\operatorname{char}(\F) \neq 2$, have been found in \cite{LM2022_2} and it turns out that 
	$$
	\operatorname{Der}(\mathfrak{d_1})=\left\lbrace \begin{pmatrix}
		x & 0 & \alpha & 0 \\
		0 & x & 0 & 0 \\
		0 & 0 & y & 0 \\
		a_1 & a_2 & a_3 & x+ y
	\end{pmatrix} \Bigg\vert x,y,\alpha,a_1,a_2,a_3 \in \mathbb{F} \right\rbrace
	$$
	Now let $D \in \operatorname{gl}(\mathfrak{d}_1)$ be an anti-derivation, then $D$ is represented, with respect to the basis $\lbrace e_1,e_2,e_3,z \rbrace$, by a matrix of the type 
	$$
	\begin{pmatrix}
		a_{11} & a_{12} & a_{13} & 0 \\
		a_{21} & a_{22} & a_{23} & 0 \\
		a_{31} & a_{32} & a_{33} & 0 \\
		A_1 & A_2 & A_3 & 0
	\end{pmatrix}
	$$
	where $D(z)=0$ because $\operatorname{Leib}(\mathfrak{d}_1)=[\mathfrak{d}_1,\mathfrak{d}_1]=\F z$, and the entries $a_{ij}$, $i,j=1,2,3$, must satisfy the following equations
	\begin{align*}
		&a_{31}+a_{32}=0,\\
		&a_{11}+a_{21}+a_{33}=0,\\
		&a_{12}+a_{22}-a_{33}=0,\\
		&-a_{12}-a_{22}-a_{33}=0.
	\end{align*}
	Thus a general anti-derivation of $\mathfrak{d}_1$ is represented by
	$$
	\begin{pmatrix}
		a_{11} & a_{12} & a_{13} & 0 \\
		a_{21} & a_{11}+a_{21}-a_{12} & a_{23} & 0 \\
		a_{31} & -a_{31} & a_{11}+a_{21} & 0 \\
		A_1 & A_2 & A_3 & 0
	\end{pmatrix}
	$$
	and, by applying the condition $(\ref{1})$, we obtain
	$$
	a_{31}=0, \; \alpha=a_{13}-a_{23}, \; y=a_{11}+a_{21}
	$$
	and
	$$
	x=a_{11}-a_{21}=a_{11}+a_{21}-2a_{12}.
	$$
	We conclude that
	$$
	\operatorname{Bider}(\mathfrak{d}_1)=\left\lbrace \left( \begin{pmatrix}
		x & 0 & \alpha & 0 \\
		0 & x & 0 & 0 \\
		0 & 0 & y & 0 \\
		a_1 & a_2 & a_3 & x+ y
	\end{pmatrix},\begin{pmatrix}
		\frac{y+x}{2} & \frac{y-x}{2} & \alpha+\beta & 0 \\
		\frac{y-x}{2} & \frac{y+x}{2} & \beta & 0 \\
		0 & 0 & y & 0 \\
		A_1 & A_2 & A_3 & 0
	\end{pmatrix} \right)  \Bigg\vert x,y,\alpha,\beta,a_i,A_j \in \mathbb{F} \right\rbrace
	$$
	and the inner biderivations are represented by the pairs of matrices of type
	$$
	\left( \begin{pmatrix}
		0 & 0 & 0 & 0 \\
		0 & 0 & 0 & 0 \\
		0 & 0 & 0 & 0 \\
		a_1 & a_1 & a_3 & 0
	\end{pmatrix},\begin{pmatrix}
		0 & 0 & 0 & 0 \\
		0 & 0 & 0 & 0 \\
		0 & 0 & 0 & 0 \\
		a_1 & -a_1 & A_3 & 0
	\end{pmatrix} \right).
	$$
\end{ex}

\section{Biderivations of low-dimensional Leibniz algebras}

Now we want to study in detail the Leibniz algebras of biderivations of low-dimensional non-Lie Leibniz algebras over a field $\F$, with $\operatorname{char}(\F) \neq 2$. There is no non-trivial Leibniz algebra in dimension 1, thus we start with two-dimensional Leibniz algebras.

\subsection{Biderivations of two-dimensional Leibniz algebras}

Let $\operatorname{dim}_{\F}L=2$, i.e.\ $L=\F e_1 + \F e_2$, then, as shown in \cite{2-dimCuvier} by C.\ Cuvier, up to isomorphism we have only two non-Lie Leibniz algebra structures on $L$.

\begin{enumerate}
	\item $L_1$ : nilpotent Leibniz algebra with non-trivial bracket $[e_2,e_2]=e_1$;
	\item $L_2$ : solvable Leibniz algebra with the table of multiplication $[e_1,e_2]=[e_2,e_2]=e_1$.
\end{enumerate}

Notice that $L_1$ is a symmetric Leibniz algebra, meanwhile $L_2$ is only a right Leibniz algebra. It turns out that
$$
\operatorname{Der}(L_1)=\left\lbrace \begin{pmatrix}
	2a & b \\
	0 & a \\
\end{pmatrix}\bigg\vert a,b \in \mathbb{F} \right\rbrace
$$
and
$$
\operatorname{Der}(L_2)=\left\lbrace \begin{pmatrix}
	a & a \\
	0 & 0 \\
\end{pmatrix}\bigg\vert a \in \mathbb{F} \right\rbrace.
$$
Moreover it is easy to check that the set of anti-derivations of $L_1$ and $L_2$ are both represented by the matrices of the form
$$
\begin{pmatrix}
	0 & x \\
	0 & y \\
\end{pmatrix}.
$$
The condition $(\ref{1})$ implies that $y=a$ for $L_1$ and $y=0$ for $L_2$, thus we have
$$
\operatorname{Bider}(L_1)=\left\lbrace \left( \begin{pmatrix}
	2a & b \\
	0 & a \\
\end{pmatrix},\begin{pmatrix}
	0 & x \\
	0 & a \\
\end{pmatrix} \right)  \bigg\vert a,b,x \in \mathbb{F} \right\rbrace
$$
and
$$
\operatorname{Bider}(L_2)=\left\lbrace \left( \begin{pmatrix}
	a & a \\
	0 & 0 \\
\end{pmatrix},\begin{pmatrix}
	0 & x \\
	0 & 0 \\
\end{pmatrix} \right)  \bigg\vert a,x \in \mathbb{F} \right\rbrace.
$$
Finally the inner biderivations of $L_1$ are represented by the pairs of matrices
$$
\left(\begin{pmatrix}
	0 & b \\
	0 & 0 \\
\end{pmatrix},
\begin{pmatrix}
	0 & -b \\
	0 & 0 \\
\end{pmatrix}\right),
$$
meanwhile the biderivations of $L_2$ are all inner.

\subsection{Biderivations of three-dimensional Leibniz algebras}

Three-dimensional complex Leibniz algebras and their derivations have been classified in \cite{3-dimLadra} and \cite{low_dim}, meanwhile the more general classification of three-dimensional right Leibniz algebras over a field $\F$, with $\operatorname{char}(\F) \neq 2$, can be found in \cite{3-dimF} and \cite{3-dimApuyov}.

Let $\dim_{\F} L=3$ and let $\lbrace e_1,e_2,e_3 \rbrace$ be a basis of $L$ over $\F$. The list of non-isomorphic three-dimensional right Leibniz algebras over $\F$ is the following.
\begin{center}
	\begin{tabular}{|l|c|}
		\hline
		\textbf{Leibniz algebra} & \textbf{Non-zero brackets} \\
		\hline
		$L_1$ & $[e_1,e_3]=-2e_1, \; [e_2,e_2]=e_1, \; [e_3,e_2]=-[e_2,e_3]=e_2$ \\
		\hline
		$L_2(\alpha), \; \alpha \neq 0$ & $[e_1,e_3]=\alpha e_1, \; [e_3,e_2]=-[e_2,e_3]=e_2$ \\
		\hline
		$L_3$ & $[e_3,e_2]=-[e_2,e_3]=e_2, \; [e_3,e_3]=-e_1$ \\
		\hline
		$L_4$ & $[e_2,e_2]=e_1, \; [e_3,e_3]=e_1$ \\
		\hline
		$L_5$ & $[e_2,e_2]=e_1, \; [e_3,e_3]=-e_1$ \\
		\hline
		$L_7(\alpha), \; \alpha \neq 0$ & $[e_2,e_2]=[e_2,e_3]=e_1, \; [e_3,e_3]=\alpha e_1$ \\
		\hline
		$L_8$ & $[e_2,e_3]=e_1$ \\
		\hline
		$L_9$ & $[e_1,e_3]=e_2, \; [e_2,e_3]=e_1$ \\
		\hline
		$L_{10}$ & $[e_1,e_3]=e_2, \; [e_2,e_3]=-e_1$ \\
		\hline
		$L_{12}(\alpha), \; \alpha \neq 0$ & $[e_1,e_3]=e_2, \; [e_2,e_3]=\alpha e_1+e_2$ \\
		\hline
		$L_{13}$ & $[e_1,e_3]=e_1, \; [e_2,e_3]=e_2$ \\
		\hline
		$L_{14}$ & $[e_1,e_3]=e_2, \; [e_3,e_3]=e_1$ \\
		\hline
		$L_{15}$ & $[e_1,e_3]=e_1+e_2, \; [e_3,e_3]=e_1$ \\
		\hline
	\end{tabular}
\end{center}
Here we use the same notation of \cite{3-dimApuyov}, but we do not report the algebras $L_6(\alpha)$ and $L_{11}(\alpha)$, where $\alpha \neq 0$, which are isomorphic to $L_4$ and $L_9$ respectively. We want to extend the results of \cite{low_dim} by completing the classification of the Lie algebras of derivations of three-dimensional Leibniz algebras over a general field $\F$, with $\operatorname{char}(\F) \neq 2$, and by finding the biderivations of this class of Leibniz algebras.

\begin{rem}
	In this section we use the following algorithm for finding derivations and anti-derivations. Let $L$ be a Leibniz algebra and let $(d,D) \in \operatorname{Bider}(L)$. Then, for every $x,y \in L$, we have
	$$
	d([x,y])=[d(x),y]+[x,d(y)], \; \; D([x,y])=[D(x),y]-[D(y),x]
	$$
	if and only if
	$$
	(d \circ \operatorname{ad}_y)(x)=(\operatorname{ad}_y \circ d)(x)+\operatorname{ad}_{d(y)}(x), \; \; (D \circ \operatorname{ad}_y)(x)=(\operatorname{ad}_y \circ D)(x)-\operatorname{Ad}_{D(y)}(x),
	$$
	thus
	$$
	[d,\operatorname{ad}_y]=\operatorname{ad}_{d(y)}, \; \; [D,\operatorname{ad}_y]=-\operatorname{Ad}_{D(y)}.
	$$
	Fixed a basis $\lbrace e_1,\cdots,e_n \rbrace$ of $L$, we have that the biderivation $(d,D)$ is represented by a pair of $n \times n$ matrices $((d_{i,j})_{i,j},(D_{i,j})_{i,j})$ and for every $i=1,\ldots,n$
	$$
	[d,\operatorname{ad}_{e_i}]=\operatorname{ad}_{d(e_i)}, \; \; [D,\operatorname{ad}_{e_i}]=-\operatorname{Ad}_{D(e_i)},
	$$
	which are equations in the entries of the matrices representing $d$ and $D$. By solving this set of equations, and after imposing the compatibility condition \textit{(1)}, we find the matrices $(d_{i,j})_{i,j},(D_{i,j})_{i,j}$.
\end{rem}

A straightforward application of the above algorithm produces the following complete classification of biderivations of three-dimensional right Leibniz algebras over $\F$. In particular we find that the dimension of these biderivation algebras lies between three and six and there are only two parameterized families of Leibniz algebra of biderivations of three-dimensional Leibniz algebras over $\F$.

\begin{thm}
	Let $\F$ be a field with $\operatorname{char}(\F) \neq 2$. The Leibniz algebras of biderivations of three-dimensional right Leibniz algebras over $\F$ can be described as follows.
	\begin{itemize}
		\item 		$\operatorname{Bider}(L_1)=\left\lbrace \left( \begin{pmatrix}
			2x & y & 0 \\
			0 & x & y \\
			0 & 0 & 0 \\
		\end{pmatrix},\begin{pmatrix}
			0 & -y & a \\
			0 & x & y \\
			0 & 0 & 0 \\
		\end{pmatrix} \right)  \Bigg\vert x,y,a \in \mathbb{F} \right\rbrace
		$;
		\item 		$\operatorname{Bider}(L_2(\alpha))=\left\lbrace \left( \begin{pmatrix}
			x & 0 & 0 \\
			0 & y & z \\
			0 & 0 & 0 \\
		\end{pmatrix},\begin{pmatrix}
			0 & a & b \\
			0 & y & z \\
			0 & 0 & 0 \\
		\end{pmatrix} \right)  \Bigg\vert x,y,z,a,b \in \mathbb{F} \right\rbrace
		$, \\ \\ if $\alpha \neq -1$ and \\ \\ $\operatorname{Bider}(L_2(-1))=\left\lbrace \left( \begin{pmatrix}
			x & 0 & 0 \\
			0 & y & z \\
			0 & 0 & 0 \\
		\end{pmatrix},\begin{pmatrix}
			0 & 0 & b \\
			0 & y & z \\
			0 & 0 & 0 \\
		\end{pmatrix} \right)  \Bigg\vert x,y,z,b \in \mathbb{F} \right\rbrace
		$;
		\item 		$\operatorname{Bider}(L_3)=\left\lbrace \left( \begin{pmatrix}
			0 & 0 & y \\
			0 & x & z \\
			0 & 0 & 0 \\
		\end{pmatrix},\begin{pmatrix}
			0 & 0 & a \\
			0 & x & z \\
			0 & 0 & 0 \\
		\end{pmatrix} \right)  \Bigg\vert x,y,z,a \in \mathbb{F} \right\rbrace
		$;
		\item 		$\operatorname{Bider}(L_4)=\left\lbrace \left( \begin{pmatrix}
			2x & y & z \\
			0 & x & 0 \\
			0 & 0 & x \\
		\end{pmatrix},\begin{pmatrix}
			0 & a & b \\
			0 & x & 0 \\
			0 & 0 & x \\
		\end{pmatrix} \right)  \Bigg\vert x,y,z,a,b \in \mathbb{F} \right\rbrace
		$;
		\item 		$\operatorname{Bider}(L_5)=\left\lbrace \left( \begin{pmatrix}
			2x & y & t \\
			0 & x & -z \\
			0 & z & x \\
		\end{pmatrix},\begin{pmatrix}
			0 & a & b \\
			0 & x & -z \\
			0 & z & x \\
		\end{pmatrix} \right)  \Bigg\vert x,y,z,t,a,b \in \mathbb{F} \right\rbrace
		$;
		\item 		$\operatorname{Bider}(L_7(\alpha))=\left\lbrace \left( \begin{pmatrix}
			\gamma x & y & z \\
			0 & x & \frac{x}{2} \\
			0 & -\frac{x}{2\alpha} & (\gamma-1)x \\
		\end{pmatrix},\begin{pmatrix}
			0 & a & b \\
			0 & x & \frac{x}{2} \\
			0 & -\frac{x}{2\alpha} & (\gamma-1)x \\
		\end{pmatrix} \right)  \Bigg\vert x,y,z,a,b \in \mathbb{F} \right\rbrace
		$, where $\gamma=\dfrac{4\alpha-1}{2\alpha}$;
		\item 		$\operatorname{Bider}(L_8)=\left\lbrace \left( \begin{pmatrix}
			x+y & z & t \\
			0 & x & 0 \\
			0 & 0 & y \\
		\end{pmatrix},\begin{pmatrix}
			0 & z & t \\
			0 & 0 & a \\
			0 & 0 & y \\
		\end{pmatrix} \right)  \Bigg\vert x,y,z,t,a \in \mathbb{F} \right\rbrace
		$;
		\item 		$\operatorname{Bider}(L_9)=\left\lbrace \left( \begin{pmatrix}
			x & y & 0 \\
			y & x & 0 \\
			0 & 0 & 0 \\
		\end{pmatrix},\begin{pmatrix}
			0 & 0 & a \\
			0 & 0 & b \\
			0 & 0 & 0 \\
		\end{pmatrix} \right)  \Bigg\vert x,y,a,b \in \mathbb{F} \right\rbrace
		$;
		\item 		$\operatorname{Bider}(L_{10})=\left\lbrace \left( \begin{pmatrix}
			x & -y & 0 \\
			y & x & 0 \\
			0 & 0 & 0 \\
		\end{pmatrix},\begin{pmatrix}
			0 & 0 & a \\
			0 & 0 & b \\
			0 & 0 & 0 \\
		\end{pmatrix} \right)  \Bigg\vert x,y,a,b \in \mathbb{F} \right\rbrace
		$;
		\item 		$\operatorname{Bider}(L_{12}(\alpha))=\left\lbrace \left( \begin{pmatrix}
			x & \alpha y & 0 \\
			y & x+y & 0 \\
			0 & 0 & 0 \\
		\end{pmatrix},\begin{pmatrix}
			0 & 0 & a \\
			0 & 0 & b \\
			0 & 0 & 0 \\
		\end{pmatrix} \right)  \Bigg\vert x,y,a,b \in \mathbb{F} \right\rbrace
		$;
		\item 		$\operatorname{Bider}(L_{13})=\left\lbrace \left( \begin{pmatrix}
			x & 0 & 0 \\
			0 & y & 0 \\
			0 & 0 & 0 \\
		\end{pmatrix},\begin{pmatrix}
			0 & 0 & a \\
			0 & 0 & b \\
			0 & 0 & 0 \\
		\end{pmatrix} \right)  \Bigg\vert x,y,a,b \in \mathbb{F} \right\rbrace
		$;
		\item 		$\operatorname{Bider}(L_{14})=\left\lbrace \left( \begin{pmatrix}
			2x & 0 & y \\
			y & 3x & z \\
			0 & 0 & x \\
		\end{pmatrix},\begin{pmatrix}
			0 & 0 & a \\
			0 & 0 & b \\
			0 & 0 & x \\
		\end{pmatrix} \right)  \Bigg\vert x,y,z,a,b \in \mathbb{F} \right\rbrace
		$;
		\item 		$\operatorname{Bider}(L_{15})=\left\lbrace \left( \begin{pmatrix}
			x & 0 & x \\
			x & 0 & y \\
			0 & 0 & 0 \\
		\end{pmatrix},\begin{pmatrix}
			0 & 0 & a \\
			0 & 0 & b \\
			0 & 0 & 0 \\
		\end{pmatrix} \right)  \Bigg\vert x,y,a,b \in \mathbb{F} \right\rbrace
		$.
	\end{itemize}
\end{thm}


\begin{thebibliography}{6}

\bibitem{loday1993version} Loday J.-L.,
"Une version non commutative des algèbres de Lie: les algèbres de Leibniz",
\emph{L’Enseignement Mathématique} (2) \textbf{39} (1993), no. 3-4, pp. 269--293.

\bibitem{blokhLie} Blokh A.,
"A generalization of the concept of a Lie algebra",
\emph{Doklady Akademii Nauk SSSR} \textbf{165} (1965), no. 3, pp. 471--473.

\bibitem{low_dim} Rakhimov I.\ and Al-Nashri A.-H.,
"On derivations of low-dimensional complex Leibniz algebras", \emph{JP Journal of Algebra, Number Theory and Applications} \textbf{21} (2011), no. 1, pp. 69--81.

\bibitem{nonabext} Liu J., Sheng Y.\ and Wang Q.,
"On non-abelian extensions of Leibniz algebras", \emph{Communications in Algebra} \textbf{46} (2018), no. 2, pp. 574--587. \url{doi: 10.1080/00927872.2017.1324870}

\bibitem{Casas} Casas J.\ M., Datuashvili T.\ and Ladra M.,
"Universal Strict General Actors and Actors in Categories of Interest", \emph{Applied Categorical Structures} \textbf{18} (2010), no. 2, pp. 85--114. \url{doi: 10.1007/s10485-008-9166-z}

\bibitem{BiggsNagy} Biggs R.\ and Nagy P.\ T.,
"On Sub-Riemannian and Riemannian Structures on the Heisenberg Groups", \emph{Journal of Dynamical and Control Systems} \textbf{22} (2016), pp. 563--594. \url{doi: 10.1007/s10883-016-9316-9}

\bibitem{ayupov2019leibniz} Ayupov S., Omirov B.\  and Rakhimov I.,
"Leibniz Algebras: Structure and Classification", \emph{CRC Press, New York} (2019), ISBN: 9781000740004.

\bibitem{erdmann2006introduction} Erdmann K.\ and Wildon M.\ J.,
"Introduction to Lie Algebras", \emph{Springer London} (2006), ISBN: 9781846284908.

\bibitem{LM2022} La Rosa G.\ and Mancini M.,
"Two-step nilpotent Leibniz algebras",
\emph{Linear algebra and its applications} \textbf{637} (2022), no. 7, pp. 119--137. \url{doi: 10.1016/j.laa.2021.12.013}

\bibitem{mac2013categories} Mac Lane S.,
"Categories for the working mathematician", Vol. 5 Graduate Texts in Mathematics, \emph{Springer-Verlag, New York} (2013), ISBN: 9781475747218.

\bibitem{LM2022_2} La Rosa G.\ and Mancini M.,
"Derivations of two-step nilpotent algebras", \emph{Communications in Algebra} (2023). \url{doi: 10.1080/00927872.2023.2222415}

\bibitem{2-dimCuvier} Cuvier C.,
"Algèbres de Leibnitz: définitions, propriétés",
\emph{Annales scientifiques de l'École Normale Supérieure} \textbf{27} (1994), no. 1, pp. 1--45.

\bibitem{3-dimLadra} Casas J.\ M., Insua M.\ A., Ladra M.\ and Ladra S.,
"An algorithm for the classification of 3-dimensional complex Leibniz algebras",
\emph{Linear algebra and its applications} \textbf{436} (2012), no. 9, pp. 3747--3756. \url{doi: 10.1016/j.laa.2011.11.039}

\bibitem{3-dimF} Rakhimov I.\ S., Rikhsiboev I.\ M. and Mohammed M.\ A.,
"An algorithm for a classification of three-dimensional Leibniz algebras over arbitrary fields", \emph{JP Journal of Algebra, Number Theory and Applications} \textbf{40} (2018), no. 2, pp. 181--198. \url{doi: 10.17654/NT040020181}

\bibitem{3-dimApuyov} Ayupov S.\ and Omirov B.,
"On 3-dimensional Leibniz algebras",
\emph{Uzbek Mathematical Journal} \textbf{1} (1999), no. 7, pp. 9--14.


\end{thebibliography}
\end{document}